\documentclass[11pt,reqno]{amsart}

\usepackage{amscd,amsmath, amsfonts, amsopn,amssymb,amsthm,multicol}
\usepackage{hyperref}
\usepackage[all]{xypic}
\usepackage{graphics}

\usepackage{hyperref}
\usepackage{here}
\usepackage{verbatim}

\def\a{\alpha}
\def\b{\beta}

\def\l{\lambda}

\def\b1{{\rm id}}

\newfont{\goth}{eufm10 scaled \magstep1}

\def\gsl{\mbox{\goth sl}}

\def\gsu{\mbox{\goth su}}

\def\gso{\mbox{\goth so}}

\def\gg{\mbox{\goth g}}
\def\gh{{\mbox{\goth h}}}
\def\gk{\mbox{\goth k}}

\newfont{\mcal}{eusm10 scaled \magstep1}

\newtheorem{Th}{Theorem}
\newtheorem{remar}[Th]{Remark}
\newtheorem{Prop}[Th]{Proposition}
\newtheorem{Cor}[Th]{Corollary}
\newtheorem{Lem}[Th]{Lemma}
\newtheorem{Def}[Th]{Definition}
\newtheorem{Ex}[Th]{Example }
\def\bt{\begin{Th}}
\def\et{\end{Th}}
\def\bp{\begin{Prop}}
\def\ep{\end{Prop}}
\def\bc{\begin{Cor}}
\def\ec{\end{Cor}}
\def\bl{\begin{Lem}}
\def\el{\end{Lem}}
\def\bd{\begin{Def}}
\def\ed{\end{Def}}
\def\bex{\begin{Ex}}
\def\eex{\end{Ex}}
\def\br{\begin{remar}}
\def\er{\end{remar}}
\def\pf{\noindent{\it Proof:\ }}
\def\qed{\hfill$\square$}

\def\be{\begin{equation}}
\def\ee{\end{equation}}
\def\ben{\begin{enumerate}}
 \def\een{\end{enumerate}}
\def\ba{\begin{array}{rlll}}
\def\ea{\end{array}}
\def\bea{\begin{eqnarray}}
\def\eea{\end{eqnarray}}
\def\bean{\begin{eqnarray*}}
\def\eean{\end{eqnarray*}}

\def\C1{cohomogeneity one  }

\begin{document}
\title{Lorentzian manifolds   with  transitive conformal  group}

\author{Dmitri Alekseevsky}

\maketitle

\centerline{A.A.Kharkevich Institute for Information Transmission
Problems}
\centerline{ B.Karetnuj  per.,19, 127051, Moscow, Russia}
\centerline{ e-mail : dalekseevsky@iitp.ru}

\begin{abstract}
 We  study   pseudo-Riemanniasn manifolds $(M,g)$ with  transitive  group  of   conformal
 transformation   which is  essential, i.e. does not preserves any metric conformal  to  $g$.
  All    such manifolds  of Lorentz   signature  with  non  exact  isotropy   representation  of  the   stability   subalgebra  are described.
   A   construction  of essential conformally   homogeneous  manifolds   with  exact  isotropy represenatation is   given. Using  spinor   formalism,  we prove    that it  provides  all
       4-dimensional  non   conformally  flat  Lorentzian    4-dimensional manifolds  with   transitive essentially   conformal  group.

\bigskip

The research is  carried out  at the IITP and is supported by an RNF grant (project n.14-50-00150).
\end{abstract}


\section{Introduction}

It is   well know   that    any  Riemannian manifold  which   admits  an  essential group of  conformal transformations
 is     conformally    equivalent  to the    standard  sphere or     the   Euclidean  space. It is   the   Lichnerowicz conjecture, proved
    in   compact  case   by  M.  Obata and J. Ferrand , and  in general  case    in  \cite{A}, \cite{A2},  \cite{Fer},\cite{F}.\\
On the other  hand, there   are many examples  of pseudo-Riemannian
	(in particular Lorentzian) manifolds   with  essential
	conformal group. Ch.\ Frances \cite{F}, \cite{F1} constructed first examples  of conformally
	essential compact Lorentzian manifolds, M.N. \	Podoksenov \cite{P}  found
	  examples  of essential conformally  homogeneous Lorentzian
	manifolds.  A  local description of   Lorentzian manifolds with
	essential  group of 	homotheties  was given in  \cite{A}.\\
	
Our   aim    is   to study  essential conformally
homogeneous pseudo-Riemannian manifolds $(M=G/H, g)$, i.e.   manifolds  with   transitive   group $G$  of  conformal transformstions
 which   does not   preserves any metric from  the   conformal  class $c=[g]$ .
We  split   all  such conformal  manifolds $(M=G/H,c)$ into two types:

A.  Manifolds   with non-exact   isotropy  representation
$$ j: \mathfrak{h} \to \mathfrak{co}(V),\, V = \mathfrak{g}/\mathfrak{h}\simeq  T_oM $$
of  the   stability   subalgebra  $\mathfrak{h}$.

B.  Manifolds   with   exact  isotropy  representation $j$.
\smallskip

We   give  a    classification of conformally   homogeneous  Lorentzian manifolds of  type A  in any  dimension
  and   classification     of non conformally  flat  manifolds  of  type  B  in   dimension  4.

 We  will  assume  that   the    transitive    conformal   group  $G$   and the  stability  subgroup  $H$   are  connected
  and   we  identify    the  pseudo-orthogonal Lie   algebra  $\mathfrak{so}_{k, \ell} = \mathfrak{so}(V)$   with   the    space $\Lambda^2V$
  of   bivectors.

\section{ Conformally    homogeneous  manifolds  and  associated  graded  Lie   algebra}
Let $(M=G/H), g)$ be  a  conformally  homogeneous
pseudo-Riemannian manifold  of  signature  $(k, \ell) =
(-\cdots-,+\cdots +)$  and $j:H\to CO(V)$ (resp., $j:\mathfrak{h}
\to \mathfrak{co}(V)$) the isotropy representation of the stability
subgroup $H$ (resp, stability subalgebra $\mathfrak{h}$) of the
point $o=eH \in M$ in the tangent space $V =T_oM$.
There is  a   filtration
$$\mathfrak{g}_{-1}= \mathfrak{g} \supset
\mathfrak{g}_{0}= \mathfrak{h}\supset \mathfrak{g}_{1} \supset
\mathfrak{g}_{2}=0$$
where $\mathfrak{g}_{1} := \mathrm{ker}j$.
The  associated transitive graded Lie algebra is

\be \bar{\mathfrak{g}}:= \mathrm{gr}(\mathfrak{g})=
\mathfrak{g}^{-1} + \mathfrak{g}^0 + \mathfrak{g}^1  =  V + \mathfrak{g}^0 +
\mathfrak{g}^1 \ee
where $V = \mathfrak{g}/\mathfrak{h}$, $\mathfrak{g}^0 :=
\mathfrak{h}/\mathfrak{g}_{1} = j(\mathfrak{h}) $ and
$\mathfrak{g}^1 = \mathfrak{g}_1 = \mathrm{ker}j $.
Transitivity  means  that
$[X,V]=0$  for   $X \in \mathfrak{g}^0 + \mathfrak{g}^{1}$ implies $X =0$.\\

\subsection{  Example:  Standard   flat model}
 Let  $\mathbb{R}^{k+1, \ell+1}$   be a pseudo-Euclidean vector  space.
 The  projectivisation
  $  S^{k, \ell} =P \mathbb{R}^{k+1, \ell+1}_0 \subset P\mathbb{R}^{k+1, \ell+1} $
   of  the isotropic   cone  $\mathbb{R}^{k+1, \ell+1}_0 \subset \mathbb{R}^{k+1, \ell+1}  $
 carries   a conformally   flat  conformal    structure $[g_{st}]$ of   signature  $(k, \ell)$. Moreover,
  $ S^{k,\ell}$  is  a   conformally homogeneous manifold  $ S^{k,\ell} = G/H=  SO_{k+1,\ell+1}/\mathrm{Sim}(V)$
  of  the pseudo-orthogonal     group $G =  SO_{k+1,\ell+1}$  of  type  A and
  the stability   subgroup  $H $ is  identified   with   the
group   of  similarities     $\mathrm{Sim}(V)=  \mathbb{R}^+\cdot SO(V)\cdot V$ of  the pseudo-Euclidean vector  space  $ V= \mathbb{R}^{k,\ell}$  via  stereographic  projection.

 The associated  graded  Lie   algebra is
\be
\mathrm{gr}(\mathfrak{so}_{k+1,\ell+1})\simeq
\mathfrak{so}_{k+1,\ell+1}= V + \mathfrak{co}(V) + V^*,
\ee
where $V^* =\mathfrak{co}(V)^{(1)}=\{T^{\xi} , [T^{\xi},X] = T^{\xi}_X = \xi(X)\mathrm{id} + X \wedge \xi   \}$
is  the  first prolongation   of   $\mathfrak{co}(V)$  and
 $   X\wedge \xi := X \otimes \xi - g^{-1}\xi \otimes gX \in \mathfrak{co}(V) $.

In the  case of Riemannian  signature $(k,\ell)=(0,n)$,
the  standard  conformal manifold  is   the  conformal sphere
$ M= S^n = SO_{1,n+1}/\mathrm{Sim}(\mathbb{R}^n)
$
\subsection{Embedding of  $gr{\mathfrak{g}}=\mathfrak{g}^{-1} + \mathfrak{g}^0 + \mathfrak{g}^1$ into  $\mathfrak{so}_{k+1, \ell+1}$}
 For   any    conformally  homogeneous  manifold  $(M= G/H, [g])$ ,  the
  associated graded Lie  algebra  $\bar{\mathfrak{g}}$ has  natural
embedding  into  the  graded Lie  algebra  $\mathfrak{so}_{k+1, \ell+1}= V + \mathfrak{so}(V)+ V^*$ as  a graded subalgebra.\\
In particular, the  conformal  structure $c$ in  $V$   induces a (may be,
degenerate) conformal structure   in  $\mathfrak{g}^1 \subset V^*$. \\

The  commutative  subalgebra  $\mathfrak{g}^1$  is   a
$\mathfrak{g}^0$-invariant  subspace
of the   first  prolongation $(\mathfrak{g}^0)^{(1)}$  and  can be
written as  $\mathfrak{g}^1 = T^{V_1^*} \subset T^{V^*}$  such
that $T^{V_1^*}\subset \mathrm{Hom}(V, \mathfrak{g}^0)$. In
particular, if $\mathfrak{g}^0 \subset \mathfrak{so}(V)$ then
$\mathfrak{g}^1=0$.

\subsection{ Subalgebras $\mathfrak{h}= \mathfrak{g}^0 \subset \mathfrak{co}(V) $   with  non  trivial  prolongation}

\bd A   decomposition
$$ V = P + E +Q $$
 of  a pseudo-Euclidean vector  space is
called  {\bf standard} if
$P,Q= P^*$  are isotropic  $k$-dimensional  subspaces  such  that $P+Q$ is  a non-degenerate
subspace  and $E$ is    the orthogonal  complement  to  $P+Q$.
\ed
We set $(P \wedge Q)^0 = \{ B \in P \wedge Q,\,   \mathrm{tr}B=0  \} = \{ \mathrm{diag}(A, -A^t),\, A \in  \mathfrak{sl}_k(\mathbb{R})\} \simeq \gsl( P)\simeq \gsl(Q)$.

\bp Let  $ \mathfrak{g}^0$  be  a proper   subalgebra    of  the
conformal linear Lie  algebra $ \mathfrak{co}(V),\,  V = \mathbb{R}^{k , \ell}$ with non-trivial
first prolongation $\mathfrak{h}^{(1)} \subset T^{V^*}$.  Then     there is  a
standard  decomposition  $V = P + E + Q$    such   that  $(\mathfrak{h})^{(1)} =
T^{g\circ P}$. \\
Moreover,  if  $k =1,\,  V = \mathbb{R}p+ E + \mathbb{R}q$, then
$$ \mathfrak{g}^0_{min}:=\mathbb{R}(\mathrm{id} - p \wedge q)+ p \wedge E\subset
\mathfrak{g}^0 \subset  \mathfrak{g}^0_{max}:= \mathfrak{g}^0_{min} + \gso(E).  $$
If $k>1$,  then
$$ \gg^0_{min}:= \mathbb{R}I+
( P\wedge Q)^0 + P \wedge (P + E) \subset \gh=\gg^0 \subset
\gg^0_{max}:=   \gg^0_{min} + \gso(E). $$
where $I= k \mathrm{id} + \mathrm{diag}(-\mathrm{id}, 0,\mathrm{id}) \in \mathfrak{gl}(P + E +Q)$.
\ep
The  proof   follows  from
\bl  If the  first prolongation  of  a     subalgebra   $\mathfrak{g}^0\subset \mathfrak{co}(V)$  contains
a non  degenerate     element  $T^{\xi},\,   g^{-1}(\xi,\xi) \neq  0$,  then
$\mathfrak{g}^0 = \mathfrak{co}(V)$.
\el
\bc  Let   $(M=G/H,c)$ be  a conformally  homogeneous  manifold.
If  the kernel $\mathfrak{g_1}$ of  the  isotropy
representation contains   a non-isotropic  element  $T^{\xi}$
 then, up  to a  covering, $M$
is conformally equivalent   to the    standard  conformal  model $(S^{k,\ell},
g_{st})$.
In  particular,  any Riemannian  conformally
homogeneous  manifold   with a non-exact isotropy  representation
is conformally   equivalent   to the  conformal sphere.
\ec

\section{Conformally homogeneous Lorentz manifolds   of  type  A}

\subsection {Conformally  flat  conformally   homogeneous  manifolds  associated   with graded subalgebra  of   $\mathfrak{so}_{k+1,\ell+1}$}

Let    $\mathfrak{g} = \mathfrak{g}^{-1} + \mathfrak{g}^0 + \mathfrak{g^1} = V + \mathfrak{g}^0 + \mathfrak{g}^1 $
be  a  graded   subalgebra  of the   graded  Lie   algebra
$$\mathfrak{so}_{k+1, \ell+1} =  V +\mathfrak{co}(V)  + V^*.$$
Assume    that $\mathfrak{g}^{1} \neq  0$ and   denote     by  $G$   the    simply connected  Lie  group associated  with   $\mathfrak{g}$  and    by  $H$  the   connected  subgroup generated   by  the subalgebra $\mathfrak{h}= \mathfrak{g}^0 + \mathfrak{g}^1$. We  assume  that  $H$ is  a  closed  subgroup.
\bt The homogeneous  manifold $M = G/H$  with   the natural     conformal  structure   defined  by  the  $j(H)$-invariant  conformal structure  in  $V$ is  a   conformally   homogeneous  manifold  of  type  A. The   commutative   subgroup  generated   by    commutative  subalgebra  $V$   has  open   orbit in $M$   and  the manifold $M$ is   conformally  flat.
\et
Note  that, in general,   the filtered   Lie  algebra $\mathfrak{g}$ of  a  conformally  homogeneous manifold is non   isomorphic   to  the  associated   graded  Lie   algebra $\bar{\mathfrak{g}}$.  In  the  next  section    we give   an  example.

\subsection{  The  standard  gradation of   $\gsu_{k+1,\ell+1}$  and  the Fefferman space}

Let $V = \mathbb{C}^{k+1,\ell+1}  = V^1 + V^0 +V^{-1} = \mathbb{C}e_+ + V^0 +
\mathbb{C}e_-$
be  a  gradation of the   complex  vector  space  $V$. We   fix  a Hermitian
form
$$V \ni Z = ue_+ + z + v e_- = (u,z,v)  \mapsto  h(Z,Z) = \bar u v +  \bar v u + h^0(z,z) $$
of   complex  signature  $(k+1,\ell+1)$   where
$h^0(z,z) =\bar{z}^t \mathbb{E}_{k,\ell}z $   is     the    Hermitian   form
in $V^0$  of  complex   signature $(k,\ell)$ with  the  Gram matrix
$\mathbb{E}_{k, \ell} = \mathrm{diag} (-1, \cdots,-1, 1,\cdots, 1)$.
   This  gradation induces   a   depth  2
gradation
of  the special unitary Lie  algebra  $\mathfrak{g}= \gsu_{k+1, \ell+1}=
\mathfrak{su}(V) = \mathfrak{aut}(V,h)$  which    may be   written   as
$$
\begin{array}{cccccc}
\mathfrak{g} & = \mathfrak{g}^{-2}& +\mathfrak{g}^{-1}& + \mathfrak{g}^0& +
\mathfrak{g}^1& + \mathfrak{g}^2 \\

\end{array}
$$

 Note  that this  gradation is
the $ \mathrm{ad}_{D}$-eigenspace   decomposition   for  $D= \mathrm{diag}(1,0,-1) =e_+ \wedge_J e_- $
where  we  use  notation  $x\wedge_Jy = x \wedge y + ix \wedge iy$.

{In matrix  notation,  the gradation  is   given    by
$$
\mathfrak{su}_{k+1,\ell+1}=
\begin{pmatrix}
\gg^0 & \gg^1 &\gg^2\\
\gg^{-1}& \gg^0& \gg^1 \\
\gg^{-2}& \gg^{-1}& \gg^0
\end{pmatrix}
=\bigg \{
\begin{pmatrix}
\l +i \mu  &  -w^* &    i\beta  \\
z &   -\frac{2i\mu}{m}\mathrm{id}+ B &   w\\
i\a & -z^* & -\l + i\mu  \\
\end{pmatrix}
\bigg \}
$$

where
$ B \in \mathfrak{su}_{k,\ell},\,  z,w  \in V^0 =  \mathbb{C}^{k,\ell} ,\, z^* := \bar z^t,\,
\alpha, \beta, \l,\mu  \in  \mathbb{R},\, m  = k + \ell.$


An   element  $L \in  \gsu_{k+1, \ell+1}$ can be  written  as
$$
\begin{array}{cccccccc}
L & = \a Q  &+ E_z &+ \mu P  &+ \lambda D& + B&+ \hat{E}_w & + \beta T  \\
\end{array}
$$
where   $D= \mathrm{diag}(1,0,-1)=  e_+\wedge_J e_- $ is  the grading element,
$$
\begin{array}{ccc}
Q&= e_-\wedge_Je_-  \in \mathfrak{g}^{-2},\,&  T=  e_+\wedge_J e_+ \in \mathfrak{g}^2\\
  E_z&= z \wedge_J e_-\in \mathfrak{g}^{-1},\,& \hat{E}_w= w \wedge_J p\in \mathfrak{g}^1\\
  P&= ip \wedge_Jq - \frac{2 i}{m} \mathrm{id_{V^0}} =i\mathrm{diag}(1, -\frac{2}{m} \mathrm{id},1) \in \mathfrak{g}^0&\\
\end{array}
$$

Denote   by $\mathcal{P} = G^0 \cdot G^+$   the  parabolic subgroup of $G =
SU_{k+1,\ell+1}$ generated  by the non-negatively graded  subalgebra
$\mathfrak{p} = \mathfrak{g}^0 + \mathfrak{g}^+= \mathfrak{g}^0 + \mathfrak{g}^1
+ \mathfrak{g}^2$.

Then  the    flag manifold  $  Fl =G/\mathcal{P} = SU_{k+1,\ell+1}/ G^0 \cdot G^+$ is the  projectivization of  the
cone of  isotropic   complex  lines  in $\mathbb{C}^{k+1,\ell+1}$.  It is
diffeomorphic  to the  sphere  if $k=0$ and it has     a natural  invariant   CR
structure. The Fefferman space  is  defined  as   the  manifold $F$ of  real
isotropic lines.  The  group   $SU_{k+1,\ell+1}$  acts transitively on
$F=SU_{k+1,\ell+1}/H $    with the  stability   subgroup $H = \mathbb{R}^+\cdot SU_{k \ell} \cdot G^+
\subset \mathcal{P}=  \mathbb{C}^* \cdot SU_{k, \ell} \cdot G^+ $.
The  Fefferman  space   is  the   total  space    of   a natural equivariant $S^1$-fibration
$$F =  SU_{k+1,\ell+1}/ H =  SU_{k+1,\ell+1}/ \mathbb{R}^+\cdot SU_{k+1. \ell+1}
\cdot G^+ \to Fl= SU_{k+1,\ell+1}/ \mathcal{P}$$
over  the  flag manifold.
The Hermitian metric $h$ of $V$ induces   an invariant   conformal metric of
signature $(2k+1, 2\ell +1))$  in $F = SU_{k+1,\ell+1}/H$.\\

The  solvable  non  commutative  Lie    algebra

$$
\mathfrak{l}  = \mathbb{R}Q + E_{\mathbb{C}^{k, \ell}}+ \mathbb{R}P=\,
\bigg \{\begin{pmatrix}
i \mu & 0 & 0\\
z& - \frac{2 \mu}{n}\mathrm{id} & 0\\
i\a &  -z^*&  i\mu
\end{pmatrix}
\bigg \}
$$
generate the  subgroup $L$ which has  an open orbit  in $F$. We  identify
$\mathfrak{l}$
with the tangent space  $T_0F$.
Then  the isotropy   representation   of  the  stability  subalgebra
$$ \mathfrak{h} =  \mathbb{R}D + \mathfrak{su_{k, \ell}} + E_{V^0} +  RT$$

is  given by
$$
\begin{array}{cll}
 j(D)   :& \a Q + E_z + \mu P &  \to  2 \a Q -E_z +0\\
  j(C) :& \a Q + E_z + \mu P & \to E_{Cz}\\
 j(\hat{E}_w):& \a Q + E_z + \mu P & \to  0+ \a E_{iw} +  \rho(w,z)P,
\end{array}
$$
where     $  C \in \mathfrak{su}_{k, \ell},\,  w^*z = \mathrm{Re}(w^*z) + \mathrm{Im}(w^*z)i = w\cdot z - \rho(w,z)i$

Note  that
$$\begin{array}{ccc}
[T,E_z]= \hat{E}_{iz}, & [T,Q]=-D, & [T,P]=0, \\
\end{array}
$$
and   that $\mathfrak{su}_{k,\ell}$ acts   by the   tautological representation  on
$E_{\mathbb{C}^{k,\ell} }$  and
$\hat{E}_{\mathbb{C}^{k,\ell}} $.
The  Fefferman space  is an example   of conformally homogeneous manifolds of  type  A,     such  that    the associated  filtered Lie algebra $\mathfrak{g}$ is not isomorphic to the   graded  Lie  algebra $\mathrm{gr}(\mathfrak{g})$. Moreover,    we have

\bt \label{theoremA}  Let  $(M=G/H,c)$ be  a  homogeneous   conformally Lorentzian
manifold of  type  A  such that the  isotropy  algebra
$j(\mathfrak{h})$ is  a proper  subalgebra of $\mathfrak{co}(V)$. If  the Lie algebra  $\mathfrak{g}$ is  not isomorphic  to  the associated    graded  Lie  algebra $\mathrm{gr}(\mathfrak{g})$,  than $M$ is  conformally isomorphic  to   the Fefferman
space  $F = SU_{1,m+1}/H$  with  conformal  metric of  signature  $(1,  2m +1)$.
\et
\subsection{ Sketch of    the proof  of    Theorem \ref{theoremA}}
\subsubsection{Step 1.}
 The graded Lie   algebra  $\mathrm{gr}(\mathfrak{g})
=\bar{\mathfrak{g}}$   associated  with $M= G/H$   has  the   form

\be \label{Lorentz graded  algebra}
\mathrm{gr}(\mathfrak{g})= \bar{\mathfrak{g}} = \bar V + (\mathbb{R}\bar D+
\bar{p} \wedge \bar E+
\bar{\mathfrak{k}} + \mathbb{R}T^{g\circ p}),\,\,\,  \bar{\mathfrak{k}} \subset  \mathfrak{so}(\bar{ E})
\ee

where $ \bar V = \mathbb{R}\bar p + \bar E+
\mathbb{R} \bar q$    is  the  standard  decompositon of   the  Minkowski vector  space
with  $ g(\bar p,\bar q)=1 $,
$\bar D:=[\bar q,T^{g \circ \bar p}]= -T^{g\circ \bar p}_{\bar q}=
-\mathrm{id}+\bar p\wedge \bar q $ .

The   element $\bar D$ defines  a  depth two gradation

$$
\begin{array}{cccccccc}
\mathrm{gr}(\mathfrak{g})= & \mathbb{R}\bar q + & \bar E +  &\mathbb{R}\bar p+
& \mathbb{R}\bar D+ & \bar{\mathfrak{k}}+
& \bar p\wedge \bar E +& \mathbb{R}T \\
\mathrm{ad}_{\bar D}  & -2            & -\mathrm{id}& 0
& 0            &   0           & \mathrm{id}          & 2
\end{array}
$$

Note that     a complementary   subspace $V$ to $\mathfrak{h}$ and a
complementary  subspace $\mathfrak{g}^0$  to $\mathfrak{h}_1$ in $
\mathfrak{h} $  defines a decomposition

\be \label{decomposition}   \mathfrak{g} = V + \mathfrak{g}^0 +
\mathfrak{h}_1 \ee

of $\mathfrak{g}$ , consistent
with the filtration $\mathfrak{g}\supset \mathfrak{g}_1
=\mathfrak{h} \supset \mathfrak{g}_1 = \mathfrak{h}_1$ and an
isomorphism of the graded vector spaces $ \mathfrak{g} $ with
$$\mathrm{gr}(\mathfrak{g}) = \bar V + \bar{\mathfrak{g}}^1 +
\bar{\mathfrak{g}}^2.$$
We will identify these spaces.

\subsubsection{Step 2}
 We can  chose  the  decomposition (\ref{decomposition}) of the
Lie
algebra $\mathfrak{g}$ such that   the  endomorphism
$\mathrm{ad}_D$ defines   a depth two gradation  as  follows

$$
\begin{array}{cccccccc}
\mathfrak{g}= & (\mathbb{R}q + & E +         &\mathbb{R}p)+ &( \mathbb{R}D+ &
\mathfrak{k}+ & p\wedge E )+& \mathbb{R}T \\
\mathrm{ad}_D  & -2            & -\mathrm{id}& 0           & 0            &
0           & \mathrm{id}          & 2
\end{array}
$$
Then  $V = \mathbb{R}q + E + \mathbb{R}p$ is a   subalgebra, which
defines a   subgroup of $G$   with open orbit.The assumptions   implies   that  $V$ is  not   commutative   subalgebra.

\subsection{Step 3}
Analyzing  Jacobi identity   we prove       that $\gg$ is of  the   following
form
\be \label{D-eigenspace decomposition}
\begin{array}{ccccccccccccccc}
\mathfrak{g}  &= &\mathbb{R}q & + &  E        & + & \mathbb{R}p &+ &
	\mathbb{R}D  &+& \mathfrak{k} & +              & p\wedge E &+ & \mathbb{R}T \\
D   & = &     -2    & + & -\mathrm{id} & + &     0      &+ &   0
	&  +&  0 &   +  & \mathrm{id} &  + & 2\\
\mathrm{ad}_p    &=&      0     & + &   A   &  + & 0  & + &  0  &   + &
	0         &+ &          A       &  + & 0 \\
\mathrm{ad}_k     &=&     0      &+ &  C  &+&     0      &+&
	0         &  +       &     \mathrm{ad}_C                   & +    &     C  & + &   0
\end{array}.
\ee
Here $k \in \mathfrak{k}$  and  $C = \mathrm{ad}_k |_{E}  \in  \mathfrak{so(E)}.$

Moreover,  $    [e,e'] = 2 \rho(e,e')q = 2 <Je,e'>q $   where  $J \in \mathfrak{so}(E)$,
 $[\mathrm{ad}_{\mathfrak{k}}, J] =0$  and   the   following  relations   hold


$$
\mathrm{ad}_T :
\begin{cases}  &  q \to  -D \\
& e \to -p\wedge e ,\, e \in E\\
& p \to 0\\
& D \to -2T\\
& \mathfrak{k} + p \wedge E \to 0,
\end{cases}
\mathrm{ad}_{p\wedge e} :
\begin{cases}                 &  q \to  -e \\
& e' \to <e,e'>p +\\
\quad &+ <Je,e'>D + K_{e,e'}\\
& p \to -p \wedge A e  \\
& D \to - p \wedge e \\
&\mathrm{ad}_{\mathfrak{k}} \ni  C \to - p \wedge \mathrm{ad}_C e , \\
& p\wedge e' \to  2<Je,e'>T.
\end{cases}
$$

Here   $K_{e.e'} \in \gk$ is  a  $\gk$-valued  symmetric   bilinear form on $E$. The
remaining  Jacobi     may be   written  as

$$	\begin{array}{cccl}
(*)  &    K_{e,e'}e''  - E_{e, e''}e'&=& -2 <Je',e''>e + <Je,e'>e'' -\\
&    & & <Je,e''>e' - <e,e'>Ae'' + <e,e''>Ae',\\
(**) &   K_{Ae,e'} + K_{e,Ae'}& =& 0,\\
(***) & C( K_{e,e'})  & = & K_{Ce,e'} + K_{e,Ce'} ,\,\,  C =
\mathrm{ad}_k,\,\,  k \in \gk .
\end{array}
$$
\subsubsection{Step 4}
The   unique  solution of  (*) is
$$    K_{e,e'}= Je \wedge e' - e \wedge Je' + <e,e'>(J-A).$$

The   equation (**) implies    that  either $J=0$ or $J $  and $A = \lambda J$    are  proportional non  degenerate  skew-symmetric   endomorphism  and   $J^2 = - p \mathrm{id}$ for  $p >0$. Rescaling metric,   we may  assume  that $J^2 = -1$.
 In  the  case $J=0$,  one  can  show  that $\mathfrak{g}$ is  isomorphic  to  $\mathrm{gr}\, \mathfrak{g} $. The   equation  (***)   shows  that $\mathrm{ad}_{ \mathfrak{k}}|_{E}  \subset \mathfrak{u}(E) = \{ C \in \mathfrak{so}(E),\,  [C,J]=0   \}  $.
Then one  can    check   that  $\gg \simeq  \gsu_{1, m+1}$   where  $n := \dim
M = 2(m+2)$    and  $M \simeq F= SU_{1,n+1}/ \mathbb{R}^+\cdot SU_n \cdot Heis(\mathbb{C}^n)$.

\subsection{The   curvature   of the Fefferman   space  and the  Cahen-Wallach   symmetric  spaces}

Recall  that   all  indecomposable Lorentzian symmetric  spaces   are  exhausted
by   the   spaces of  constant  curvature
  and the  Cahen-Wallach  symmetric  spaces $CW^{1,n-1}_S$
Let
$$V= \mathbb{R}^{1,n-1} =  \mathbb{R}q + E + \mathbb{R}p$$
be the  standard  decomposition of   the Minkowski   space  and   $e_i$    an  orthonormal basis  of $E$.  Then the
contravariant  curvature  tensor  $R_S$ of  the  Cahen-Wallach space  is  given by
$$  R_S = \sum_{i=1}^{n-2} q\wedge S e_i \vee q \wedge e_i .$$
It   defines     a Lie    algebra  with  a  symmetric  decomposition
$$   \mathfrak{g} = \mathfrak{h}+ V = q \wedge E + V \subset \gso(V) +V $$
with  the   Lie  bracket     $[x,y] = - R(x,y) \in \mathfrak{h}  =  q \wedge  E ,\, x,y \in V$.
The  Cahen-Wallach  space $CW^{1,n-1}_S  =  G/H$ the symmetric   manifold  associated  with  this  symmetric  decomposition.
  It is
  conformally  flat if  and only if
$S = \lambda \mathrm{id}$, see \cite{G}.

\bt   For  any point  $x$  of  the Fefferman   space   $(F, [g])$  there is    a metric   $g \in [g]$
whose     contravariant curvature tensor  at   $x$ coincides  with the     curvature tensor  of the conformally  flat
Cahen-Wallach   space.   In  particular, the  Fefferman   space  is   conformally  flat.
\et

\section{ Petrov  classification  of  Weyl   tensors }
\subsection{Spinor formalism}
To  describe 4-dimensional Lorentzian  conformally
homogeneous manifolds of type  B, we recall a spinor
description of
Weyl tensor of  a  Lorentzian  4-manifold.\\

Let $\mathbb{S} $ be the complex 2-space with the symplectic form
$\omega = e_- \wedge e_+$  where    $e_+,e_-$  is   a symplectic basis   of   $\mathbb{S}$ and  we identify   $\mathbb{S}$ with  the   dual  space $\mathbb{S}^*$.
$\omega(e_+ , e_-)=1$    which   is  identified    with  the  dual  space
  The associated standard basis $E_- = E_{21}, \, E_0 = E_{11}- E_{22}, E_+ =
E_{12}$ of  the unimodular Lie  algebra  $\mathfrak{sl}_2(C)$
 defines   a gradation
$$\mathfrak{sl}_2(C)=\mathfrak{g}^{-1}+ \mathfrak{g}^0 + \mathfrak{g}^1=
\mathbb{C} E_- +
\mathbb{C}(E_0) + \mathbb{C}E_+. $$
The   space $\mathbb{S}\otimes \bar{ \mathbb{S}}$ of  Hermitian
bilinear  forms  has   the  basis
$e_i \otimes \bar{ e_j},\, i,j \in \{ +,- \}$
where  $\bar e_+, \bar e_-$ is  the basis of  the  conjugated vector
space $\bar {\mathbb{S}} =\bar{ \mathbb{C}^2} = \{\bar{ z_{+}}{\bar
	e_+} + \bar {z_-} \bar {e_-} \}$. If  $j : a \otimes \bar b \mapsto
(a\otimes \bar b)^* = b\otimes \bar a$ is   the  Hermitian
conjugation, then the fix point  space  $V = (\mathbb{S} \otimes
\bar {\mathbb{S}})^j$ of $j$ is  the  space of  Hermitian  symmetric
matrices.

We   may write
$$
\begin{array}{cc}
V =&  \{ X = uE_{1\bar 1}+ (z E_{1\bar 2} + \bar z E_{2 \bar 1}) + v
E_{2 \bar 2}\}  = \big \{
\begin{pmatrix}
u&z  \\
\bar z & v
\end{pmatrix}, u,v \in \mathbb{R}, z \in \mathbb{C}\big \}.
\end{array}
$$

We  set $p =2E_{1\bar 1}, \, q = 2 E_{2 \bar 2}, E_z = z E_{1\bar 2}
+ \bar z E_{2 \bar 1}$  such that $E = \{E_z,\, z \in \mathbb{C}\}
\simeq \mathbb{C}$ and denote by
$$  V = V^{-1}+ V^0 + V^1 = \mathbb{R}q + E + \mathbb{R}p  $$
the associated  gradation  of $V$. The  determinant
defines  the  Lorentz  metric   in $V$ :
$$g(X,X)=  X\cdot X  =  \det X   =  uv  - z \bar z  = uv - x^2 -y^2,\,  z = x
+i y$$
such that
$$p^2:= p \cdot p = q^2 =0,\, p \cdot q =2, \, e_1^2 = e_i^2 =-1, e_1
\cdot e_i =0,
(\mathbb{R}p + \mathbb{R}q ) \perp E$$
where $e_1 :=E_1, e_i = E_i$. \\
For $X,Y \in V $  we  denote  by  $X \wedge Y: Z \mapsto <Y,Z>X -
<X,Z>Y$  the  associated  endomorphism  from $\mathfrak{so}(V)$.
The   group  $SL(\mathbb{S})$  acts isometrically   in    $V$  by
$$\varphi : SL(\mathbb{S}) \ni A  \mapsto \phi(A):  X \mapsto
AXA^*.$$
The  associated  isomorphism of  Lie  algebras
$\mathfrak{sl}(\mathbb{S})$  and  $\mathfrak{so}(V)$ is   given by
$$
\begin{array}{cc}
\varphi(E_0) = 2p\wedge q & \varphi(iE_0)=2e_1 \wedge e_i, \\
\varphi(E_+) = \sqrt{2}e_1\wedge p& \varphi(iE_+)=-\sqrt{2}e_i \wedge p, \\
\varphi(E_-) = \sqrt{2}e_1\wedge q & \varphi(iE_-)=- \sqrt{2} e_i \wedge q.
\end{array}
$$

\subsection{Spinor  description of  the space  $\mathcal{R}_0(V)$  of Weyl
	tensors}

Recall that  the  space of Weyl tensors is  defined  by
$$\mathcal{R}_0(V) =\{W \in \mathrm{Hom}(\Lambda^2V, \mathfrak{so}(V)),
\mathrm{cycl}\,W(X\wedge Y)Z =0,$$
$$\mathrm{tr}(X \to W(X,\cdot)X)=0, \forall X,Y,Z \in V   \}.$$
Recall that
$\Lambda^2 V \simeq
\mathfrak{sl}_2(\mathbb{C}) \simeq \mathbb{C}^3$
where   the  complex  structure in $\Lambda^2V$  is defined  by Hodge  star operator.
Note  that
$  V^{\mathbb{C}}= \mathbb{S}\otimes \bar{ \mathbb{S}}$
and
$ \Lambda^2V^{\mathbb{C}}  =  S^2\mathbb{S}\otimes \bar{\omega}  + \omega
\otimes {S^2}(\bar{\mathbb{S}})$
where $\omega, \bar{ \omega}$ are  symplectic forms in $\mathbb{S}$
and $\bar{ \mathbb{S}}$.\\
  We    denote by
$S^2_0(\Lambda^2(V))$   the 5-dimensional complex space  of  trace
free symmetric  complex  endomorphisms  of  the  complex   space
$\Lambda^2(V) = \mathbb{C}^3$.

\bt (A. Petrov, R. Penrose)   There  exists   a  natural
isomorphisms of  $\mathfrak{sl}_2(\mathbb{C})$- modules
$$  \mathcal{R}_0(V) \simeq S^2_0(\Lambda^2(V)) = S^4 (\mathbb{S}^*).$$
The covariant form  $ g\circ W$  of  the  Weyl tensor $W$
associated  with   symmetric  4-form  $\varphi  $  is  given by
$$ W_{\varphi} = \varphi \otimes \bar{\omega}^2 +  \omega^2 \otimes \bar{
	\varphi}. $$
\et
\subsection{Petrov  classification of  Weyl tensors}

Since  any symmetric   form $\phi \in S^4(\mathbb{S})$  can be
factorized  into a  product of  linear  form
$   \phi =  \alpha \beta \gamma \delta $
we  get  the  following  classification of Weyl  tensors:\\
Type (4) or (N) $ \phi = \alpha^4 $; \,\,  Type (31) or (III) $ \phi = \alpha^3 \beta $;\\
Type (22) or (D) $ \phi = \alpha^2 \beta^2$;\,\, Type  (211) or (II) $ \phi = \alpha^2 \beta \gamma$ ;\\
Type  (1111) or (I)  $\phi = \alpha \beta \gamma \delta$, \\
where $\alpha, \beta, \gamma, \delta$ are  different linear forms in
$\mathbb{S}$.
Each linear form  $\alpha$ in  spinor  space $\mathbb{S}$ up to a
scaling is  defined  by  its kernel $\alpha =0$
which is  a point  in to projective  line $\mathbb{C}P^1 = S^2$. So
up to a  complex  factor,  the 4-form $\phi$  is  determined  by  4
points  on the  conformal  sphere.
For a  symmetric 4-form $\phi $   we  denote  by
$\mathfrak{aut}(\phi)$ (respectively,$\mathfrak{conf}(\phi) )$ the Lie
subalgebra of
$\mathfrak{sl}_2(\mathbb{C})$  which preserves $\phi$ ( respectively, preserves
$\phi$ up to a  complex
factor).
\bp \label{Weylspinors}

i)  $\mathfrak{conf}(\phi)=0$  for a  form of types (1111), (211);\\
ii) $\mathfrak{aut}(\phi)=0$  for a  form of  types  different from
$(2,2)$  and $(4)$;\\
iii)   $\mathfrak{conf}(\phi)= \mathbb{C} =
\mathbb{C}E_0$   for a form of types $(31)$;\\
iv)  $\mathfrak{conf}(\phi)  = \mathfrak{aut}(\phi) = \mathbb{C}E_0$ for
type $(22)$ ;\\
v) $\mathfrak{conf}(\phi) = \mathbb{C}E_0 +
\mathbb{C}E_+,\,   \mathfrak{aut}(\phi) = \mathbb{C}E_+$  for  type
$(4)$.\\
In particular, only the  form  of  type $(31)$  and  $(4)$  admits a  conformal   transformation  which  is not  an  authomorphism.
\ep

\pf There   exists  unique   conformal transformation of  the sphere
which transform  three   different  points into   another  three
different points. This implies  the  first claim. If $\phi =
\alpha^4$,  then   the  stabilizer of $\phi$ in
$\mathfrak{sl}_2(\mathbb{C})$  is  the  same  as  the  stabilizer of
the 1-form $\alpha$. We  may assume that  $\alpha =e_-^*= (0,1)$.
Then $ \mathfrak{aut}(\phi) = \mathbb{C}E_+$   and
$\mathfrak{conf}(\phi) = \mathbb{C}E_0 + \mathbb{C}E_+$. If $\phi =
\alpha^2 \beta^2 $ or $\alpha^3 \beta$, we  may assume that $\alpha,
\beta$ are basic 1-forms and  then   the stabilizer of
$\mathbb{C}\phi$ will be the diagonal  subalgebra. In the  first
case it   preserves $\phi$. \qed

\section{Conformally  homogeneous manifolds of  type B}

In  this  section  we describe a class   of    conformally   homogeneous
pseudo-Riemannian manifolds   of  type  B  and    prove that  all   4-dimensional
    conformally   homogeneous  non   conformally  flat manifold  belong   to  this  class.

\bp  Let   $M=G/H$   be  a  conformally   homogeneous manifold  of  type  B.
 Then   the  isotropy Lie algebra    $j(\mathfrak{h}) \subset  \mathfrak{co}(V),\,  V = T_0M$  has  a  decomposition

 $$j(\mathfrak{h})  = \mathbb{R}D +
\mathfrak{l},$$
where $\mathfrak{l} \subset \mathfrak{so}(V)$
is  an
ideal  of  $\mathfrak{h\mathbb{}}$ an  the  endomorphism  $D = \mathrm{id} +
C,\, C \in \mathfrak{so}(V)$ is  a non  trivial  homothety.
  \ep

\pf  Indeed,assume  that  $j(\mathfrak{h})\subset  \mathfrak{so}(V)$. Then   the isotropy group   $j(H)$ preserves   a  metric $g_0$ in  the   tangent  space   $V = T_0M$ which   can be  extended   by left translations   to   $G$-invariant metric $g$  on the    homogeneous    space   $M = G/H$. Hence,  the conformal group   $G$ is not  essential.
\qed

\subsection{ A   construction  of pseudo-Riemannian   conformally   homogeneous manifold of  type  B }

Let $V = \mathbb{R}q + E + \mathbb{R}p$     be  a    standard  decomposition of  a
pseudo-Euclidean  vector  space $(V, g = <.,.>)$  of   signature $(k,\ell)$.
The   homothety
$D = \mathrm{id} + q \wedge p   \in   \mathfrak{co}(V)$.
 defines  a  gradation $V =
\mathbb{R}p + E + \mathbb{R}q=  V^0 + V^1 + V^2$. A non-degenerate
2-form $\omega(x,y)$   in  $E$     defines
the    structure of    the  Heisenberg
Lie   algebra    with the  center   $\mathbb{R}q$ and the bracket  $[x,y] = \omega(x,y)q,\,  x,y\in  E$
in  $ \mathfrak{heis}(E) =  E + \mathbb{R}q$ . Moreover,   an   endomorphism  $A
\in \mathrm{End}(E)$   with
$$(A\cdot \omega)(x,y) := \omega(Ax,y) + \omega(x,Ay) = \lambda  \omega(x,y)$$
is a  derivation of  this  algebra   and  defines   the  structure  of   a
graded  Lie  algebra
$$ V =V^0 + V^1 +V^2 = \mathrm{R}p   + 	\mathfrak{heis}(E),$$
such  that  $\mathrm{ad}_p q=\lambda q, \,
\mathrm{ad}_p|_E =A$,
with the  grading   element    $D = \mathrm{id} + q \wedge p$.    Denote
by  $G$  the Lie  group  generated  by the  Lie     algebra
$\gg  = \mathbb{R}D + V $ and  by  $H$   the  closed  subgroup  generated  by the
subalgebra  $\mathbb{R}D$.
\bp    The  metric $g$  in   $V$ defines   an invariant  pseudo-Riemannian   conformal  structure in the manifold   $M = M(\lambda, \omega, A)=  G/H$.    The manifold   $M$  is a  conformally    homogeneous
manifold  of  type B.\\
  The curvature  operator     of  the manifold   $M$   is   given  by
  $$R_{pq}= R_{qx}=0,\,
   R_{px}=  (A^aA^s -  A^sA - JA^s)x \wedge q,\,\,  x\in E $$

   where  $g^{-1}\circ \omega  = 2 J$       and  $A^a =\frac12(A+A^t), \, A^s =\frac12(A-A^t) $
    are skew-symmetric and symmetric  parts of  $A$.
   In particular,   in  general   the manifold  $M$ is    not  conformally flat.
   \ep

\subsection{Classification  of Lorentzian   4-dimensional conformally   homogeneous  manifolds of  type  B}

\bt \label{theoremB} Any  conformally homogeneous  4-dimensional Lorentzian  manifold of  type  B which is
not   conformally  flat  is   conformally isometric  to   a manifold
$M(\lambda,\omega, A)$.
\et

The  proof   is based on

\bl If     a conformally homogeneous Lorentzian 4-manifold   $M^4=G/H$  of  type  B  is  not    conformally  flat,    then  the  isotropy  Lie  algebra
 contains the homothety  $D  = \mathrm{id}+ q\wedge p$  with respect  to   an   approprite   standard  decomposition
$V = \mathbb{R}p + E + \mathbb{R}q $
 of  the  tangent  space    $V = T_o M  $.
\el

\pf  Let   $D = \mathrm{id} + C  \in  j(\mathfrak{h})$  be a non  trivial homothetic  endomorphism,
$C \in \mathfrak{so}(V)$.
 By  assumption,the  Weyl tensor  $W  \neq 0$. Since $\mathrm{id}  \cdot W =-2$   and
$D \cdot W = (\mathrm{id} + C) \cdot W = 0,\,  C\cdot W = 2 W. $
Then    $C\cdot \phi=2\phi$ ,  where
  $ \phi \in  S^4(\mathbb{S}^2)$  is  4-form   which  represents  $W$ .
  Proposition \ref{Weylspinors}  shows  that   the  4-form $\phi$
has    Pertov  type (4)  or  $(31)$  and     $C = -\frac{1}{2}E_0 + b E_-  \in \mathfrak{sl}_2(\mathbb{C})$. Changing the basis,  we  may   assume  that  $b=0$.  Then     $\varphi^{-1}(C) = q \wedge  p \in  \mathfrak{so}(V)$   and    $D = \mathrm{id} + q\wedge p$.
\qed

\end{document}